\documentclass{conm-p-l}

\usepackage{setspace}

\theoremstyle{definition}

\theoremstyle{remark}

\numberwithin{equation}{section}

\newcommand{\bee}{\begin{equation}}
\newcommand{\pa}{\partial}
\newcommand{\al}{\alpha}

\begin{document}

\title{The Legacy of the IST}

\author{David J. Kaup}
\address{Department of Mathematics, University of Central Florida, Orlando, FL 32816-1364}
\curraddr{}
\email{kaup@ucf.edu}
\thanks{The author was supported in part by NSF Grant \#0129714. The author thanks an anonymous referee for his comments, and also H. Steudel for his comments.}


\date{}



\begin{abstract}
We provide a brief review of some of the major research results arising from the method of the Inverse Scattering Transform.
\end{abstract}

\maketitle









\section{Introduction}

I will give a brief review of several items in the Legacy of the Inverse Scattering Transform. In no way is this to be a complete review, since the Legacy has become so vast. However, I will treat those items with which I am most familiar, and try to detail their significance and importance.

There is no doubt that the most important contribution was the famous classical Gardner, Greene, Kruskal and Miura (GGKM) work  \cite{GGKM67} of 1967 on the KdV equation. This was the starting point. They had found a very strange and new method for solving the initial value problem of a nonlinear evolution equation, the KdV. At that time, and even for several years later, this strange new method was considered to be only a novelty, since it would only work for that one equation, the KdV. Shortly thereafter, as a prelude to what was to follow, Peter Lax  \cite{L68}  showed that if given an appropriate linear operator, $L$, dependent on a potential, $u(x)$, then one could always construct an infinite sequence of evolution operators, $B$, each of which would satisfy
\bee\label{e:one}
BL-LB=\pa_t L.
\end{equation}
This sequence of evolution operators could be generated by simply increasing the order of the spatial differentials contained in $B$. Then from (1.1) one would obtain additional nonlinear evolution equations, each of the form
\bee
\pa_t u=K(u)
\end{equation}
where $K$ was some (nonlinear) operator. All these additional higher order evolution equations would be solvable by this same technique. This collection is now known as the KdV hierarchy.

If we consider the eigenvalue problem for $L$, 
\begin{subequations}
\bee
L\psi =-\lambda \psi ,
\end{equation}
where $\lambda$ is the eigenvalue, and append to it the condition
\bee
\pa_t \psi =B\psi ,
\end{equation}
\end{subequations}
then it is easy to see that (1.1) is simply the integrability condition for (1.3). Furthermore, as Lax pointed out for the KdV case, (1.1)-(1.3) also implies that the eigenvalues, $\lambda$, in (1.3a) would be stationary,
\bee
\pa_t\lambda =0,
\end{equation}
a relation that would occur time-and-time again as the study of integrable equations would expand in the decades to follow.

It was also about this time that the term ``radiation" was introduced. We haven't said anything yet about solitons or solitary waves, but more will be said later. For now, let us note that a remarkable feature of the GGKM method of solution was the appearance of fully nonlinear solitary wave solutions, called solitons. The other part of the solution has been called ``radiation", and is essentially linear-like in its behavior. The asymptotics (long-time behavior) of the total solution are generally that the radiation does disperse away, leaving the solitons traveling in a sea of decaying radiation.

As to nomenclature, we shall refer to $(L+\lambda )\psi$ as the eigenvalue problem, $\psi$ as the eigenfunctions, and $B$ as the evolution operation. The pair $[L+\lambda ,B]$ is known as the ``Lax pair". For the KdV equation, the Lax pair is
\begin{subequations}
\bee
(L+\lambda )\psi =(\pa_x^2 +u+\lambda )\psi =0
\end{equation}
\bee
\pa_t\psi =b\psi =(\al -4\pa_x^3 -6u\pa_x -3u_x )\psi
\end{equation}
where $\alpha$ is an arbitrary constant and the integrability condition is the KdV equation:
\bee
\pa_tu+\pa_x^3 u+6u\pa_xu=0 .
\end{equation}
\end{subequations}

We note here, given $L$ and $B$, it follows that one can then obtain $K(u)$. However, an important problem is given $K(u)$, construct $L$ and $B$. The solution of this inverse problem is still an area of active research. One method that sometimes works for this is called ``Painlev\'e Analysis". For a description of this aspect of the Legacy, the reader is referred to Choudhury's article in this same issue.

It was not until 1971, that the next physically significant integrable system was uncovered by V. Zakharov and A.B. Shabat (ZS) \cite{ZS71}, which was the focusing Nonlinear Schr\"odinger Equation (NLS)
\bee
iq_t=-q_{xx}-2q^* q^2 .
\end{equation}
This equation required a different eigenvalue problem,
\begin{subequations}
\begin{align}
v_{1x}+i\zeta v_1 &= qv_2\\
v_{2x}-i\zeta v_2 &= rv_1
\end{align}
\end{subequations}
where $(v_1,v_2)$ is the eigenvector, $\zeta$ is the eigenfunction, and $q$ and $r$ are the ``potentials". For the focusing NLS case, one has $r=-q^*$ and $r=+q^*$ for the defocusing case. The time evolution operator, $B$, is given by \cite{ZS71}:
\begin{subequations}
\begin{align}
i\pa_tv_1 &= -i(4\zeta^2+2qr)v_1+(4\zeta q+2iq_x)v_2\\
i\pa_t v_2 &= (4\zeta r-2ir_x)v_1+i(4\zeta^2+2qr)v_2 .
\end{align}
\end{subequations}
Exactly as was shown by Lax \cite{L68} for the KdV, one also has a hierarchy here, which can be obtained by generalizing (1.8) to higher orders in $\zeta$. In 1972, Wadati \cite{W72} found the next member of this hierarchy, the ``modified KdV" (mKdV)
\bee
q_t+q_{xxx}+6q^2q_x=0
\end{equation}
which was also integrable. Its eigenvalue problem was again the ZS eigenvalue problem, (1.7), but where now $r=-q$, and $q$ real. Also, (1.8) had to be generalized to be cubic in $\zeta$.

By this time, it was becoming apparent to many researchers, that this strange method found by GGKM was not simply a novelty. Rather, there was something very significant underlying all of this. This became even more obvious when Ablowitz, Kaup, Newell and Segur (AKNS) presented a method of solution of both the Goursat and Cauchy initial value problems of the sine-Gordon equation \cite{AKNS73a}. This was also based on the ZS eigenvalue problem, but with a very different form for the $B$ in (1.8): it was now inversely proportional to the spectral parameter,  $\zeta$. The sine-Gordon equation was well known at that time. It had a long history, first occuring in 1853 in differential geometry, and was the first equation for which B\"acklund tranformations and N-solitons solutions were found. It was known in solid state physics in the 1930's, and in 1965 had found applications in optics. 

The IST solution of the sine-Gordon equation was shortly followed by another letter \cite{AKNS73b} pointing out how one could generate a large number of integrable equations, each of which were physically significant and important. With one general approach, AKNS were able to reproduce all the Lax pairs found up to that time, and were able to connect the form of the dependence of $B$ on $\zeta$ to the linear dispersion relation, $\omega (k)$.  (The linear dispersion relation relates how the frequency, $\omega$, depends on the wave vector, $k$, in the linear limit, where plane waves, $e^{i(kx-\omega t)}$, are the natural solutions.)  In 1974, they published their classic AKNS paper \cite{AKNS74}, wherein they described in detail this new method of solution, calling it the method of the Inverse Scattering Transform (IST). One of the major points of this classic was that the IST could be viewed as a nonlinear extension of the method of the Fourier Transform. 

This was also the start of the explosion in research on solitons and integrable systems, because unbeknownst to most westerners, Faddeev, Zakharov and their students were all very busy in the same direction. In the next few years, many important papers were to be published on the IST and related issues.

\section{The Legacy}

Beginning in 1974, it becomes difficult to try to detail all the results. Nevertheless, we will now discuss in general terms, the legacy which followed from this. In the following, we will list the general areas of the legacy, and briefly describe the importance and the major contributions made to each one.
 
\subsection{Method of Solution --- the IST}

Above all, the IST is a method of solution for integrable nonlinear equations. It was the pioneering work of GGKM~\cite{GGKM67}, ZS \cite{ZS71} and AKNS \cite{AKNS74} which made the most significant impact and set the tone which followed. Consider what the IST allows one to do. One may take any reasonable initial data, and by means of the direct scattering transform (the eigenvalue problem), transform the initial data into ``scattering data". For the KdV, the scattering data consists of a reflection coefficient $[\rho(k);-\infty <k<\infty ]$, bound state eigenvalues $[\kappa_j ;j=1,2, \dots ,N]$ and bound state normalization coefficients $[C_j ;j=1,2, \dots ,N]$ where $N$ is the total number of bound states, usually finite. Now, Lax \cite{L68} showed that the eigenvalues, $\kappa_j$, would be independent of time, due to (1.1). GGKM showed that if $u(x,t)$ evolved according to the KdV equation, the reflection coefficient, $\rho$, and the normalization coefficients, $C_j$, would evolve according to
\begin{subequations}
\begin{align}
\pa_t \rho (k;t) &= \rho (k;0)e^{8ik^3t} , \\
\pa_t C_j(t) &= C_j(0)e^{-8\kappa_j^3t} .
\end{align}
\end{subequations}
Thus it becomes a very simple matter to determine the scattering data at any later time.

Next, one used the solution of the inverse scattering problem, which is the core of this method of solution, to reconstruct the potential(s). One transforms (with the IST) from the scattering data, at time $t$, back to the potential(s), at time $t$. For the KdV, the necessary steps are to first construct
\bee
F(z;t)=\frac{1}{2\pi} \int_{-\infty}^\infty \rho (k;t)e^{ikz}dk+\sum_{j=1}^N C_j(t) e^{-\kappa_j z} ,
\end{equation}
then one solves the linear integral equation
\bee
K(x,y)+F(x+y)+\int_x^\infty K(x,s)F(s+y)ds=0,
\end{equation}
for $K(x,y;t)$. Lastly $u(x,t)$ is constructed from
\bee
u(x,t)=-2\frac{dK(x,x;t)}{dx} .
\end{equation}
All integrable systems are solved by an IST of the above format, although there can be a wide variation in the form of the formats. Some eigenvalue problems are higher order and/or even multidimensional. But there is always some scattering problem which maps the potential(s) into a set of scattering data (the Direct Scattering Transform). There is always some evolution of the scattering data, as in (2.1). There is always an Inverse Scattering Transform that allows one to map from the scattering data back to the potential(s), as in (2.2)-(2.4).

One may describe this method of solution as a method for solving the initial value problem of a nonlinear equation by using only linear techniques. Furthermore, one may also say that one solves these nonlinear problems by \underline{not} solving the nonlinear problem. Instead one solves two related \underline{linear problems}. 

To understand these comments, consider
\bee
L(u, \zeta )\cdot V=0,
\end{equation}
where $L$ is the eigenvalue problem in the Lax pair, as in (1.1)-(1.3) and $\zeta$ is a spectral parameter. Let the second component of the Lax pair be $B(u,u_x, \dots ,\zeta ,\pa_x , \dots )$, where
\bee
B\cdot V=\pa_tV.
\end{equation}
Now note that (2.5) determines the $x$-dependence of $V(x,t;\zeta )$, whereas (2.6) determines the $t$-dependence. Thus one function, $V(x,t;\zeta )$, is being determined by two equations. In general, such would overdetermine $V$, and in order for a mutual solution to exist, certain consistency conditions must be satisfied. This condition is the single integrability condition, (1.2). When (1.2) is satisfied, then a mutual solution exists for (2.3)-(2.6).

But if this integrability condition is nothing more than the nonlinear evolution equation, it therefore follows that if, by some means, we can construct a single-valued solution for $V(x,t;\zeta )$, for some $u(x,t)$, which satisfies each component of the Lax pair, then it follows that $u(x,t)$ must satisfy the nonlinear evolution equation. So, one could say that the entirety of the method of the IST is based on \underline{not} solving integrable nonlinear equations (at least, not directly). Instead, we solve them indirectly, exactly by the same format used in any transform method. This is as follows.  We satisfy (2.5) by mapping $u(x,t)$ into scattering data and constructing the eigenfunctions. We can do this for certain classes of potentials (i.e. $L_1 \bigcap L_2$). Then we always will have reflection coefficient(s), as a function of the spectral parameter, $\zeta$, and certain bound state data (eigenvalues and normalization coefficients).  This is a linear problem. We satisfy (2.6) by requiring the scattering data to evolve appropriately (as in (2.1) for the KdV). This is usually easy to do, since scattering data is typically defined for $x\rightarrow \pm \infty$, where potentials normally approach specified values (usually zero) in that limit. So one really only needs $B(x\rightarrow \pm \infty$).  This is also a linear problem. Therefore, by fixing the scattering data to evolve appropriately, we have effectively forced  $u(x,t)$ to evolve by the nonlinear evolution equation. To solve for $u(x,t)$, we need to solve the inverse scattering problem,  which is also a linear problem. The construction of the kernel(s), $F(z)$, as in (2.2), is another linear problem. Equation (2.3) is a linear integral equation for $K(x,y)$. Then $u$ is reconstructed as in (2.4) by a linear operation. Thus with only linear techniques, we are able to solve these nonlinear evolution equations.

\subsection{Soliton Solutions}

One of the unique features of the IST is that it allows one to construct an infinity of exact nonlinear solutions, called the $N$-soliton solutions. These are also called reflectionless potentials because the scattering data consists of only bound state scattering data, with all reflection coefficients set equal to zero. In the case of the KdV equation, when the reflection coefficients vanish, the function $F(x+y;t)$ in Eq. (2.3), then separates into a finite sum of products of known functions of $x$ and $y$, allowing one to obtain $K(x,y)$, also as a finite sum of known functions. The same is true for all other cases integrable by the IST.

The value of these solutions is tremendous, since they allow one to study and obtain exact results for these systems. The most important solution in this class is always the one-soliton solution, since it is the basic building block of these solutions, and also of any interactions between solitons and between solitons and  radiation. The 1-soliton solution of the KdV hierarchy is
\bee
u=\frac{2\eta^2}{\text{cosh}^2 \{ \eta [x-x_0(t)]\}}
\end{equation}
and that of the NLS hierarchy is
\bee
q=\frac{2\eta e^{-2i\xi [x-x_0(t)]}e^{-i\phi_0(t)}}{\text{cosh} \{2\eta [x-x_0(t)]\}}
\end{equation}
where the form of $x_0(t)$ and $\phi_0(t)$ depends on the member of the hierarchy.
The 2- and 3-solitons solutions demonstrate soliton interactions and collisions. Typically what happens in any soliton collision is that after the collision, the solitons separate out according to their individual velocities, completely unscathed, except for a possible phase shift in their positions and/or phase. There are two well known exceptions to this. In the soliton decay case of the 3WRI \cite{KRB79}, an initial soliton in the high frequency wave can decay into its two daughter waves, transferring its identity to them. In the vector NLS \cite{M73}, similar inelastic type collisions occur whereby a soliton in one mode (polarization), can flip to another combination of modes. Bound states can occur in the sine-Gordon field, where one can have stable 2-soliton bound states, called ``breathers" \cite{AKNS73a}, which are localized oscillations of the sine-Gordon field.

\subsection{Hierarchies}

Given any eigenvalue problem, Lax \cite{L68} had noted that one could always take the eigenvalue problem (the first component of the Lax Pair), and by simply extending the order of the evolution operator, $B$, one could generate another nonlinear integrable equation. All these nonlinear integrable equations which have a common eigenvalue problem in the Lax Pair, but different evolution operators, $B$, is called a ``hierarchy". Thus, there is a hierarchy for every one of these eigenvalue problems. For the Schr\"odinger equation, the most important members are the KdV equation, a 5th order KdV equation \cite{K80} and one-dimensional ``caviton" equation \cite{K88} (the analogy of the sine-Gordon equation for the Schr\"odinger eigenvalue problem). In the ZS problem, we have a ``workhorse" as far as physical equations are concerned. If we allow $r$ in (1.7) to be in general independent  from $q$, then in addition to the NLS \cite{ZS71}, the hierarchy contains the modified KdV \cite{W72}, the sine-Gordon \cite{AKNS73a}, the sinh-Gordon equation \cite{AKNS73a}, coherent pulse propagation and self-induced transparency (SIT) \cite{AKN74}, stimulated Raman scattering (SRS) \cite{CS75}, and the defocusing NLS \cite{ZS73}. The hierarchy containing the three-wave resonant interaction (3WRI) includes all three forms of this interaction (explosive, soliton decay, and stimulated backscatter (SBS)) \cite{KRB79,ZM73,K76a}, as well as the Manakov vector NLS~\cite{M73}.

\subsection {Inverse Scattering}

Another aspect of the legacy is the wide variety of inverse scattering problems solved. In 1967, the IST of the Schr\"odinger equation had just recently been obtained \cite{GL51,FS63}, and it was only in 1972 that the IST of the ZS eigenvalue problem had been solved \cite{ZS71}. Since that time, there has been a multitude of other and even more complex scattering problems solved.

The next one was the solution of the third-order eigenvalue problem for the three wave resonant interaction (3WRI) \cite{ZM73,K76a}, which interaction we shall return to later. This one differed from the ZS significantly only in the additional order of the problem. An important subcase of this was the inverse scattering solution for the eigenvalue problem for the vector nonlinear Schr\"odinger (VNLS) equation, solved by Manakov \cite{M73}. This VNLS equation is a very important and key equation for several studies of nonlinear optical pulses propagating in optical fibers \cite{LK97, HWK97}. The generalization of the order 3 problem to order $n$ has been done by Gerdjikov and Kulish \cite{GK81}.

More complex forms have also been done. The first one of these was the eigenvalue problem for the sine-Gordon equation in laboratory coordinates \cite{FT74, K75}. Here one has the spectral parameter distributed among various potential terms, and is
\begin{subequations}
\begin{align}
v_{1x}+ \left( \frac{i}{2} \zeta -\frac{i}{8\zeta} \cos u \right) v_1 &= \left[ \frac{i}{8\zeta} \sin u- \frac14 (u_x+u_t) \right] v_2 \\
v_{2x}- \left( \frac{i}{2} \zeta -\frac{i}{8\zeta} \cos u \right) v_1 &= \left[ \frac{i}{8\zeta} \sin u- \frac14 (u_x+u_t) \right] v_1 . 
\end{align}
\end{subequations}
One notes that due to this form, one cannot easily express this equation in the standard form as $L\cdot V=\lambda V$ where $L$ is a nondegenerate differential operator, $V$ is the eigenvector and $\lambda$ is the eigenvalue. However, in spite of this, one still can solve the direct and inverse scattering problems for such systems. Further examples of such eigenvalue problems are the eigenvalue problems for the massive Thirring model \cite{KuMi77}, and the derivative NLSL \cite{KN78b}. A more standard form is the cubic generalization of the Schr\"odinger equation \cite{K80, Zak74}, 
\bee
\psi_{xxx} +6Q\psi_x +6R\psi=\lambda \psi
\end{equation}
which appeared as the eigenvalue problem for the Boussinesq, Sawada-Kotera equation and the Kaup-Kuperschmidt equation. There is also the inverse scattering for multidimensional problems, such as the 3D form of the 3WRI \cite{RK81}, the KPI and II equations \cite{FA83a}, as well as the DSI and II equations \cite{FS90}. Several aspects of these are still of current research interest.

\subsection{Perturbations and Closure}

Once one has an exact method for the solution of a system, it then becomes possible to develop perturbation methods, to study nearby systems. This work was first done in 1976 for the ZS eigenvalue problem for the one-soliton solution \cite{K76b}, with a general summary of the perturbation method being given in 1978 \cite{KN78a}.

Key to this, is the concept called ``closure" or ``completeness", which itself arises from the one-to-one nature of the direct scattering transform and the IST~\cite{K76c}. What this simply means is that given any potential in the appropriate class, there exist a unique set of scattering data that can be associated with it by the eigenvalue problem of the Lax pair, and vice versa for the IST. Thus for any potential whose evolution is slightly perturbed away from its integrable value, by the direct scattering transform, it will be mapped into some other scattering data near the initial integrable scattering data. The perturbation problem is then to determine how this  scattering data in the perturbed case  evolves in time. Once its evolution is known, then by the IST, one may map back to  the potentials and then obtain their evolution.

One solves this by relating variations in the potentials to variations in the scattering data, with the transformation from the former to the latter being accomplished by the so-called ``squared eigenfunctions". Then with these squared eigenfunctions, one may obtain the evolution of the scattering data under a given perturbation.

Under perturbations, one no longer has the simple evolution of the scattering data as in (2.1). Rather, one has a slow mixing of the various elements of the scattering data: solitons will decay and/or be pumped, transforming some of their energy into radiation and/or absorbing energy from the perturbations. In addition, radiation modes will similarly grow and/or decay. Many aspects of this have been covered in the review of soliton perturbations done by Kivshar and Malomed \cite{KM89} as well as in Ref. \cite{KN78a}.

However, whenever the eigenvalue problem has a singular structure, one has difficulties. The classical example is the perturbation theory of the KdV, where the Schr\"odinger eigenvalue problem is singular at $k=0$. What happens is that for a pure one-soliton solution, $\rho (k)=0$ at $k=0$. But if there is just even the smallest amount of radiation present, we have $\rho (k=0)=-1$. Now, the width of this region in $k$-space can be very small (proportional to the radiation density), so this region could be vanishingly small and one might have to look very close to $k=0$ in order to see it.  However, due to this singular behavior, any perturbations of the KdV soliton which will create any radiation, will always have secular terms~\cite{KM77}. Physically what is occurring is that a shelf or a depression is forming due to the perturbation \cite{KN78a}, and this action does  generate a \underline{finite} shift in the scattering data as the perturbation vanishes.

\subsection{General Integrable Evolution Equations}

Given a linear dispersion relation, what are the possible integrable evolution equations for that system? Well, for any integrable system solvable by the ZS IST or the Schr\"odinger IST, that answer has been given by Newell and Kaup \cite{KN79}. They showed that the most general system is any of the AKNS polynomial forms \cite{AKNS74}, coupled with a generalized SIT system. Using the properties of the ZS squared eigenstates and their closure, they were able to construct the most general evolution equation, given the linear dispersion relation. In other words, the dispersion relation of the linear theory determines the nonlinear terms, given the eigenvalue problem. The same could be done for any other hierarchy.

Now, a burning question has always been ``Is it ever possible for an exactly integrable system to have evolving eigenvalues (i.e. for soliton amplitudes to evolve)?" In general, the answer to this is ``No". However, Kaup and Newell did find exceptions. Such equations can indeed be constructed, but potential applications for them seem to be remote.

\subsection{Optical Systems}

One question which arose naturally after a few years is ``Why are so many of the physical integrable systems related to nonlinear optical systems?" As one goes down the list, one has the sine-Gordon equation (which applies to two-level atoms and is the sharp line limit of SIT), nonlinear Schr\"odinger (focusing and defocusing --- the NLS is almost always the weakly nonlinear limit of any almost monochromatic envelope \cite{N80}), self-induced transparency (SIT --- also more generally referred to as two-level coherent propagation), three-level coherent propagation \cite{M84}, three-wave resonant interactions (3WRI) \cite{ZM73,K76a,KRB79}, second harmonic generation (SHG) \cite{K78}, the three-dimensional form of the 3WRI \cite{ZM76,C79,K81}, stimulated Raman scattering (SRS) \cite{CS75,K83,FM99}, two photon propagation (TPP) \cite{S72,S83,S88}, and degenerate TPP (DTPP) \cite{SK96}. This is not meant to be a comprehensive list, but it does include the major integrable nonlinear optical systems.

Of these, the 3WRI was the first system to demonstrate a major departure from the accepted and expected soliton behavior. The first deviation was that the radiation would never asymptotically vanish. Also after any and all collisions, nor would the radiation separate out from the solitons.  The reason for this is that the 3WRI has \underline{no} dispersion, and therefore solitons will never separate from the radiation present. They each have the same velocity. One could now ask what is the significance of the 3WRI solitons, when they differ so much from KdV or NLS solitons? To answer that the best we can say is that in the 3WRI, any solitons seem to simply represent a ``packet" of something. It is a unit which cannot be broken up. Although there are certain exchange rules for the exchanging of solitons between the three envelopes \cite{K76a,KRB79}, nevertheless solitons in the 3WRI do seem to be some robust and coherent part of the envelope. On the other hand, the radiation component of any envelope has no such finite ``packet" size, but rather can always be subdivided and redistributed among the three envelopes, subject only to the Conservation of Action laws \cite{KRB79}.

It was the SRS system that first brought to the forefront certain ambiguities with the IST on finite intervals. The  first solution of the finite interval case, by using an infinite interval IST, was given in \cite{S83}, where the general IST for SRS was developed, and features of the solution were discussed. Numerics of SRS have been studied by Hilfer and Menyuk \cite{HM90}, the asymptotical form of the solution was described by Kaup \cite{K92}, and was later solved as a Riemann-Hilbert problem by Fokas and Menyuk  \cite{FM99}. In the meantime, by analytically solving a model initial value problem, Menyuk and Kaup essentially found that for all these integrable nonlinear optical problems  (except NLS and 3WRI), one could just as easily describe the solution as ``an infinity of solitons with no radiation", as well as by the usual description of ``a finite number of solitons in a sea of radiation". Briefly, why this could happen is basically the same reason as why a Fourier transform on a finite interval has a variety of forms. There one could use either a cosine series, or a sine series, or the exponential series to represent the function. For a given function, the coefficients in each of these series is quite different. Another way to look at this, is that on a finite interval, there is an infinite number of ways to take and combine plane waves to reconstruct some function \underline{inside the finite interval}. Similarly for the IST, on a finite interval, there is no unique form for the scattering data.

Pursuing this further, consider the case where solitons are forbidden, as in the ZS $r=+q^*$ case. Now, what is going to happen in this case where no solitons are allowed on the infinite interval, if we try to represent the solution on the finite interval with no reflection coefficients, and only with solitons? What happens to this system on a finite interval? Well, we again find something surprising. Taking SHG as an example case of these defocusing systems, we find that the solution of a simple initial value problem can also be given in terms of no radiation, but now (since regular solitons are forbidden)  an infinity of what is called ``virtual solitons"~\cite{KS01}. These are indeed interesting objects. 

The possibility of their existence was noted way back in 1974 \cite{AKNS74}.  However, no known application was then known for them. In the ZS IST, the typical soliton for the focusing case $(r=-q^*)$ is of the form
\bee
q=\frac{2 \eta e^{-2i\xi (x-x_0)}e^{-i\phi_0}}{\text{cosh} [2\eta (x-x_0)]} \; .
\end{equation}
Now, if one simply assumes that the defocusing case $(r=+q^*)$ has one bound state in the scattering data, one obtains the ``virtual soliton" solution
\bee
q=\frac{2\eta }{\text{sinh} [2\eta (x-x_0)]} \; ,
\end{equation}
which clearly will always be singular on any infinite interval. However, on a finite or a semi-infinite interval, as long as $x_0$ lies outside the physical region, the solution is perfectly valid. Thus a virtual soliton is the singular, ZS, $r=+q^*$, soliton. (For certain technical reasons, they also were first found in the lower half complex plane~\cite{KS01}, which is another reason for the prefix ``virtual". For more details on this, see Steudel's contribution in this same series.)

As to my first question as to why so many nonlinear optical systems are integrable, perhaps the principle reason for this is the vast orders of magnitude difference between the speed of light and other velocities in these physical systems (which are typically acoustic velocities, and/or electronic drift velocities). Due to this disparity in velocities, one can then expect a multiple-scale expansion to give quite good results. Also one would then expect to see the higher order terms scaling like some power of this velocity ratio, and therefore rapidly vanishing.

\subsection{Benjamin-Ono Equation}

Although we could perhaps have included the Benjamin-Ono (BO) equation in the section on eigenvalue problems, it is unique enough to justify some additional comments. First it is a one dimension eigenvalue problem but it is also a \underline{nonlocal} eigenvalue problem \cite{FA83b}. Actually it can be formulated as an electrostatics problem (Poisson's equation) in two dimensions,  since it can be stated in terms of Hilbert transforms and their analytical properties. This is probably the simplest explanation as to why Fokas and Ablowitz \cite{FA83b} could term the BO equation as a ``pivotal equation for multidimensional problems". It indeed does contain this multidimensional flavor. It also should be noted that there are now two versions of its IST \cite{FA83b, KM98}. This is not entirely surprising, since the multidimensional 3D-3WRI \cite{K81} also contains a multitude of different forms of the IST. The multidimensional flavor of the BO equation stands out even more when one considers its closure property, its perturbation theory, and its Hamiltonian structure \cite{KLM98}.

\subsection{Reduction Problems}

It seems that almost any integrable system can be found in the Yang-Mills field, if one knows how to find the right reduction~\cite{YM}. However, there is another quite useful direction that one can take in reducing integrable systems. Let us first note that if one simplifies, or reduces the number of degrees of freedom of  an integrable system by some set of constraints, consistent with the integrable flows, then the reduced system will also be integrable.  As one example, consider the propagation of $N$ solitons in an optical fiber. It is very important to maintain the spacings of these $N$ solitons over long distances. So a key question is the stability of such an arrangement, which is simply an $N$-soliton solution. However, to try to analyze a 100-soliton state, where every soliton has approximately the same amplitude, is not an easy thing to analyze, since the transmission coefficient, $a$, has a zero of order 100 at the pole. Another approach for $N$ larger than 2 or 3 is clearly needed. One way to do this is to approximate the system as a lattice, since in general, one expects the solitons to have an almost equally spacing (or being absent if representing a zero) and of equal amplitude. Then it turns out that when the $N$-soliton system is reduced to a lattice system, it reduces to a Complex Toda Lattice of $N$ points, which is integrable. One now can study the stability problem of the reduced system, the Complex Toda Lattice, and transfer the results to the $N$-soliton optical pulse \cite{GEKDU98}. From this, one obtains general criteria about how the phase of each successive soliton should be adjusted, to maximize the stability of the pulse train, and etc.

\subsection{The Fokas Method}

Another new aspect of the IST has been recently developed by Prof. Fokas \cite{F01}. There are complex practical problems requiring the solution of 2-dimensional, linear, constant coefficient, PDEs but with complex boundaries, such as wedges and polygons. Before his work, there was no general analytical technique for constructing solutions with such complicated boundaries. However, it is a simple matter to create a Lax pair for almost any linear PDE with constant coefficients \cite{BF02}. Once this is done, one then can approach the solution of these linear two-dimensional, boundary value problems from the IST point of view, whereby one solves (or satisfies) both components of the Lax pair, thereby also satisfying the integrability condition, which is just the PDE to be solved. The method can also be extended to integrable nonlinear PDEs and evolution equations, however the solution then frequently requires the solution of a Riemann-Hilbert problem for the reflection coefficients.

\subsection{Unsolved Problems}

Now we will discuss what the hard problems are. They are hard because they haven't been solved. This is not a complete list, but is the start of a list of problems in need of solutions.

Although we hae a solution for the IST of the cubic eigenvalue problem, (2.20), there still is no equivalent of the GLM equations for this system. What we do have is the solution of the Riemann-Hilbert problem for the eigenfunctions. What is next needed is a representation of the eigenfunctions in terms of some transformation kernels (like the $K(x,y;t)$ for the KdV). Once these are known, then the equivalent GLM equations will follow. Some progress in this direction has recently been made by A. Parker \cite{P01}.

There are still aspects of KPI and II that are of interest. Numerical simulations~\cite{ISS95} provide valuable insight into the evolution of these equations. I would also say that there are probably still some unanswered questions about DSI and II. However, the unfortunate thing about these two equations is that potential applications seem to be almost lacking, due to the scales involved. (Typical parameters for DS solitons in water require meter-like distances horizontally, but only centimeter-like water depths, and even smaller wave amplitudes.)

Questions about perturbation theory and closure relations seem to have become quite well understood, based on the works of Gerdjikov, Ivanov and Kulish \cite{GIK80}, Beals and Coifman \cite{BC84}, and J. Yang \cite{Y02}. They are also quite well understood for the 3D3WRI and DSI, since they both use the same eigenvalue problem. However, I know of no published work in this area.

Of the optical problems there is still interest in SHG and DTPP, and the latter  has even the direct scattering problem to be detailed, as well as the inverse scattering problem.

I also strongly suspect that there are other integrable systems still to be found. Here one would have to apply Painlev\'e analysis and see what will result.



\begin{thebibliography}{99}
\bibitem{GGKM67} C.S. Gardner, J.M. Greene, M.D. Kruskal and R.M. Miura, Method for solving the Korteweg-deVries equation, {\it Phys. Rev. Lett.} {\bf 19} (1967), 1095-1097.
\bibitem{L68} P.D. Lax, Integrals of nonlinear equations of evolution and solitary waves, {\it Comm. Pure Appl. Math.} {\bf 21} (1968), 467-490.
\bibitem{ZS71} V.E. Zakharov and A.B. Shabat, Exact theory of two-dimensional self-focusing and one-dimensional self-modulation of waves in nonlinear media, {\it Zh. Eksp. Teor. Fiz} {\bf 61}, (1971), 118 [{\it Sov. Phys. JETP} {\bf 34}, (1972), 62].
\bibitem{W72} M. Wadati, The exact solution of the modified Korteweg-de Vries equation, {\it J. Phys. Soc. Jap.} {\bf 32} (1972), 1681.
\bibitem{AKNS73a} M.J. Ablowitz, D.J. Kaup, A.C. Newell, and H. Segur, Method for solving the sine-Gordon equation, {\it Phys. Rev. Letters}, {\bf 30}, 25 (1973), 1262-1264.
\bibitem{AKNS73b} M.J. Ablowitz, D.J. Kaup, A.C. Newell, and H. Segur, Nonlinear-evolution equations of physical significance, {\it Phys. Rev. Letters} {\bf 31}, 2 (1973), 125-127.
\bibitem{AKNS74} M.J. Ablowitz, D.J. Kaup, A.C. Newell and H. Segur, The inverse scattering transform-Fourier analysis for nonlinear problems, {\it Studies in Appl. Math.} {\bf LIII}, 4 (1974), 249-315.
\bibitem{KRB79} D.J. Kaup, A. Reiman, and A. Bers, Space-time evolution of nonlinear three-wave interactions. I. Interaction in a homogeneous medium, {\it Rev. Mod. Phys.} {\bf 51} (1979), 275-310.
\bibitem{M73} S.V. Manakov, On the theory of two-dimensional stationary self-focusing of electromagnetic waves, {\it Zh. Eksp. Teor. Fiz} {\bf 65} (1973), 505-516 [{\it Sov. Phys.-JETP} {\bf 38}, 2 (1974), 248-253].
\bibitem{K80} D.J. Kaup, On the inverse scattering problem for cubic eigenvalue problems of the class $\psi_{xxx}+6Q\psi_x+6R_\psi =\lambda\psi$, {\it Studies in Appl. Math.} {\bf 62} (1980), 189-216.
\bibitem{K88} D.J. Kaup, Cavitons are solitons: An integrable ponderomotive system, {\it Phys. Fluids} {\bf 31}, 6 (1988), 1465-1470.
\bibitem{AKN74} M.J. Ablowitz, D.J. Kaup, and A.C. Newell, Coherent pulse propagation, a dispersive, irreversible phenomenon, {\it J. Math. Phys.} {\bf 15}, 11 (1974), 1852-1858.
\bibitem{CS75} F.Y.F. Chu and A.C. Scott, Inverse scattering transform for wave-wave scattering, {\it Phys. Rev. A} {\bf 12} (1975), 2060.
\bibitem{ZS73} V.E. Zakharov and A.B. Shabat, Interaction between solitons in a stable medium, {\it Zh. Eksp. Teor. Fiz} {\bf 64} (1973), 1627-1639.
\bibitem{ZM73} V.E. Zakharov and S.V. Manakov, Resonant interaction of wave packets in nonlinear media, {\it Zh. Eksp.  Teor. Fiz. Pis'ma Red.} {\bf 18}, (1973), 413 [{\it Sov. Phys.-JETP Lett.} {\bf 18} (1973), 243].
\bibitem{K76a} D.J. Kaup, The three-wave interaction --- A nondispersive phenomenon, {\it Studies in Appl. Math.} {\bf 55} (1976), 9-44.
\bibitem{GL51} I.M. Gel'fand and B.M. Levitan, On the determination of a different equation from its spectral function, {\it Izvest. Akad Nauk S.S.S.R. ser. matem} {\bf 15} (1951), 309-360 [{\it Amer. Math. Soc. Transl.}, Ser. 2, 1 (1955), 259-309].
\bibitem{FS63} L.D. Faddeyev and B. Seckler, The inverse problem in quantum theory of scattering, {\it J. Math. Phys.} {\bf 4} 1 (1963), 72-104.
\bibitem{LK97} T.I. Lakoba and D.J. Kaup, Perturbation theory for the Manakov soliton and its applications to pulse propagation in randomly birefringent fibers, {\it Phys. Rev. E} {\bf 56}, 5 (1997), 6147-6165.
\bibitem{HWK97} H.A. Haus, W.S. Wong, and F.I. Khatri, Continuum generation by perturbation of soliton, {\it J. Opt. Soc. Am. B} {\bf 14}, 2 (1997), 304-313.
\bibitem{GK81} V.S. Gerdjikov and P.P. Kulish, The generating operator for the $n\times n$ linear system, {\it Physica} {\bf 3D} (1981), 549-564.
\bibitem{FT74} L.D. Faddeev and L.A. Takhtajan, Essentially nonlinear one-dimensional model of the classical field theory, {\it Teor. Mat. Fiz.} {\bf 21} (1974), 160-174.
\bibitem{K75} D.J. Kaup, Method for solving the sine-Gordon equation in laboratory coordinates, {\it Studies in Appl. Math.} {\bf 54}, 2 (1975), 165-179.
\bibitem{KuMi77} E.A. Kuzhetsov and A.V. Mikhailov, On the complete integrability of the two-dimensional classical Thirring model, {\it Teor. Mat. Fiz.} {\bf 30}, 3 (1977), 303-314.
\bibitem{KN78b} D.J. Kaup and A.C.Newell, An exact solution for a derivative nonlinear Schr\"odinger equation, {\it J. Math. Phys.} {\bf 19} (1978), 798-801.
\bibitem{Zak74} V. Zakharov, On the stochastization of one-dimensional chains of nonlinear oscillators, {\it Sov. Phys. JETP} {\bf 38} (1974), 108-110. 
\bibitem{RK81} A. Reiman and D.J. Kaup, Multi-shock solutions of random phase three-wave interactions, {\it Physica} {\bf 3D} (1981), 389-394.
\bibitem{FA83a} A.S. Fokas and M.J. Ablowitz, On the inverse scattering of the time-dependent Schr\"odinger equation and the associated Kadomtsev-Petviashvili (I) equation, {\it Stud. Appl. Math} {\bf 69} (1983), 211-228.
\bibitem{FS90} A.S. Fokas and P.M. Santini, Dromions and a boundary value problem for the Davey-Stewartson 1 equation, {\it Physica D} {\bf 44} (1990), 99-130.
\bibitem{K76b} D.J. Kaup, A perturbation expansion for the Zakharov-Shabat inverse scattering transform, {\it SIAM J. Appl. Math.} {\bf 31}, 1 (1976), 121-133.
\bibitem{KN78a} D.J. Kaup and A.C. Newell, Solitons as particles, oscillators, and in slowly changing media: A singular perturbation theory, {\it Proc. R. Soc. Lond.} A {\bf 361} (1978), 413-446.
\bibitem{K76c} D.J. Kaup, Closure of the squared Zakharov-Shabat eigenstates, {\it J. Math. Anal. Appl.} {\bf 54}, 3 (1976), 849-864.
\bibitem{KM89} Y.S. Kivshar and B.A. Malomed, Dynamics of solitons in nearly integrable systems, {\it Rev. Mod. Phys.} {\bf 61}, 4 (1989), 763-915.
\bibitem{KM77} V.I. Karpman and E.M. Maslov, Perturbation theory for solitons, {\it Zh. Eksp. Teor. Fiz.} {\bf 73} (1977), 537-559.
\bibitem{KN79} D.J. Kaup and A.C. Newell, Evolution equations, singular dispersion relations, and moving eigenvalues, {\it Adv. Math.} {\bf 31}, 1 (1979), 67-100.
\bibitem{N80} A.C. Newell, The inverse scattering transform, in {\it Solitons}, Ed. R.K. Bullough and P.J. Caudrey, pp. 176-242 [Springer-Verlag, New York (1980)].
\bibitem{M84} A.I. Maimistov, Rigorous theory of self-induced transparency in the case of a double resonance in a three-level medium, {\it Sov. J. Quantum Electron.} {\bf 14}, 3 (1984), 385-389.



\bibitem{K78} D.J. Kaup, Simple harmonic generation: An exact method of solution, {\it Studies in Appl. Math.} {\bf 59} (1978), 25-35.
\bibitem{ZM76} V.E. Zakharov and S.V. Manakov, The theory of resonance interaction of wave packets in nonlinear media, {\it Zh. Eksp. Teor. Fiz.} {\bf 69} (1975), 1654 [{\it Sov. Phys.-JETP} {\bf 42} (1976), 842].
\bibitem{C79} H. Cornille, Solutions of the nonlinear 3-wave equations in three spatial dimensions, {\it J. Math. Phys.} {\bf 20} (1979), 1653.
\bibitem{K81} D.J. Kaup, The solution of the general initial value problem for full three dimensional three-wave resonant interaction, {\it Physica D} {\bf 3} (1981), 374-395.
\bibitem{K83} D.J. Kaup, The method of solution for stimulated Raman scattering and two-photon propagation, {\it Physica D} {\bf 6} (1983), 143-154.
\bibitem{FM99} A.S. Fokas and C.R. Menyuk, Integrability and self-similarity in transient stimulated Raman scattering, {\it J. Nonlinear Sci.} {\bf 9} (1999), 1-31.
\bibitem{S72} H. Steudel, Stimulated Raman scattering with ultrashort pulses, {\it Exp. Tech. Phys.} {\bf 20} 5 (1972), 409-415
\bibitem{S83} H. Steudel, Solitons in stimulated Raman scattering and resonant two-photon propagation, {\it Physica 6} {\bf D} (1983), 155-178.
\bibitem{S88} H. Steudel, $N$-soliton solutions to degenerate self-induced transparency, {\it J. Mod. Opt.} {\bf 35}, 4 (1988), 693-702.
\bibitem{SK96} H. Steudel and D.J. Kaup, Degenerate two-photon propagation and the oscillating two-stream instability: The general solution for amplitude-modulated pulses, {\it J. Mod. Opt.} {\bf 43}, 9 (1996), 1851-1866.
\bibitem{HM90}G. Hilfer and C.R. Menyuk, Stimulated Raman scattering in the transient limit, {\it J. Opt. Soc. Am.} B {\bf 7} (1990), 739.
\bibitem{K92} D.J. Kaup, The asymptotic solution of the stimulated Raman-scattering equation, in {\it Painlev\'e Transcendents}, Ed. D. Levi and P.Winternitz, Plenum Press, New York, 1992, 345-351.
\bibitem{KS01} D.J. Kaup and H. Steudel, Virtual solitons and the asymptotics of second harmonic generation, {\it Inv. Prob.} {\bf 17} (2001), 959-970.
\bibitem{FA83b} A.S. Fokas and M.J. Ablowitz, The inverse scattering transform for the Benjamin-Ono equation~--- A pivot to multidimensional problems, {\it Stud. Appl. Math} {\bf 68} (1983), 1-10.
\bibitem{KM98} D.J. Kaup and Y. Matsuno, The inverse scattering transform for the Benjamin-Ono equation, {\it Stud. Appl. Math.} {\bf 101} (1998), 73-98.
\bibitem{KLM98} D.J. Kaup, T.I. Lakoba and Y. Matsuno, Complete integrability of the Benjamin-Ono equation by means of action-angle variables, {\it Phys. Letters A} {\bf 238} (1998), 123-133.
\bibitem{YM} M.J. Ablowitz, S. Chakravarty and R. Halburd, The generalized Chazy equation from the self-duality equations, Preprint (2002).
\bibitem{GEKDU98} V.S. Gerdjikov, E.G. Evstatiev, D.J. Kaup, G.L. Diankov and I.M. Uzunov, Stability and quasi-equidistant propagation of NLS soliton trains, {\it Phys. Letters A} {\bf 241} (1998), 323-328.
\bibitem{F01} A.S. Fokas, Two-dimensional linear PDEs in a convex polygon, {\it Proc. R. Soc. Lond. A} (2001).
\bibitem{BF02} D. ben-Avraham and A.S. Fokas, The solution of the modified Helmholtz equation in a triangular domain and an application to diffusion-limited coalescence (Preprint, 2001)
\bibitem{P01} A. Parker, A reformulation of the ``dressing method" for the Sawada-Kotera equation, {\it Inverse Problems}, submitted.
\bibitem{ISS95} E. Infeld, A. Senatorski and A.A. Skorupski, Numerical simulations of Kadomtsev-Petviashvili soliton interactions, {\it Phys. Rev. E} {\bf 51} (1995), 3183-3191.
\bibitem{GIK80} V.S. Gerdzhikov, M.I. Ivanov and P.P. Kulish, Quadratic bundle and nonlinear equations, {\it Teor. Mat. Fiz.} {\bf 44}, 3 (1980), 342-357.
\bibitem{BC84} R. Beals and R.R. Coifman, Scattering and inverse scattering for first order systems, {\it Comm. Pure Appl. Math.} {\bf 37} (1984), 38-90.
\bibitem{Y02} J. Yang, Private communication (2002).
\end{thebibliography}
\end{document}